\input amstex
\input amsppt.sty

\magnification1200
\hsize13cm
\vsize19cm

\catcode`\@=11
\font@\twelverm=cmr10 scaled\magstep1
\font@\twelveit=cmti10 scaled\magstep1
\font@\twelvebf=cmbx10 scaled\magstep1
\font@\twelvei=cmmi10 scaled\magstep1
\font@\twelvesy=cmsy10 scaled\magstep1
\font@\twelveex=cmex10 scaled\magstep1

\newtoks\twelvepoint@
\def\twelvepoint{\normalbaselineskip15\p@
 \abovedisplayskip15\p@ plus3.6\p@ minus10.8\p@
 \belowdisplayskip\abovedisplayskip
 \abovedisplayshortskip\z@ plus3.6\p@
 \belowdisplayshortskip8.4\p@ plus3.6\p@ minus4.8\p@
 \textonlyfont@\rm\twelverm \textonlyfont@\it\twelveit
 \textonlyfont@\sl\twelvesl \textonlyfont@\bf\twelvebf
 \textonlyfont@\smc\twelvesmc \textonlyfont@\tt\twelvett
%Erg"nzung des fetten Small-Capitals-Fonts:
%
 \ifsyntax@ \def\big##1{{\hbox{$\left##1\right.$}}}%
  \let\Big\big \let\bigg\big \let\Bigg\big
 \else
  \textfont\z@=\twelverm  \scriptfont\z@=\tenrm  \scriptscriptfont\z@=\sevenrm
  \textfont\@ne=\twelvei  \scriptfont\@ne=\teni  \scriptscriptfont\@ne=\seveni
  \textfont\tw@=\twelvesy \scriptfont\tw@=\tensy \scriptscriptfont\tw@=\sevensy
  \textfont\thr@@=\twelveex \scriptfont\thr@@=\tenex
        \scriptscriptfont\thr@@=\tenex
  \textfont\itfam=\twelveit \scriptfont\itfam=\tenit
        \scriptscriptfont\itfam=\tenit
  \textfont\bffam=\twelvebf \scriptfont\bffam=\tenbf
        \scriptscriptfont\bffam=\sevenbf
  \setbox\strutbox\hbox{\vrule height10.2\p@ depth4.2\p@ width\z@}%
  \setbox\strutbox@\hbox{\lower.6\normallineskiplimit\vbox{%
        \kern-\normallineskiplimit\copy\strutbox}}%
 \setbox\z@\vbox{\hbox{$($}\kern\z@}\bigsize@=1.4\ht\z@
 \fi
 \normalbaselines\rm\ex@.2326ex\jot3.6\ex@\the\twelvepoint@}

\font@\fourteenrm=cmr10 scaled\magstep2
\font@\fourteenit=cmti10 scaled\magstep2
\font@\fourteensl=cmsl10 scaled\magstep2
\font@\fourteensmc=cmcsc10 scaled\magstep2
\font@\fourteentt=cmtt10 scaled\magstep2
\font@\fourteenbf=cmbx10 scaled\magstep2
\font@\fourteeni=cmmi10 scaled\magstep2
\font@\fourteensy=cmsy10 scaled\magstep2
\font@\fourteenex=cmex10 scaled\magstep2
\font@\fourteenmsa=msam10 scaled\magstep2
\font@\fourteeneufm=eufm10 scaled\magstep2
\font@\fourteenmsb=msbm10 scaled\magstep2
\newtoks\fourteenpoint@
\def\fourteenpoint{\normalbaselineskip15\p@
 \abovedisplayskip18\p@ plus4.3\p@ minus12.9\p@
 \belowdisplayskip\abovedisplayskip
 \abovedisplayshortskip\z@ plus4.3\p@
 \belowdisplayshortskip10.1\p@ plus4.3\p@ minus5.8\p@
 \textonlyfont@\rm\fourteenrm \textonlyfont@\it\fourteenit
 \textonlyfont@\sl\fourteensl \textonlyfont@\bf\fourteenbf
 \textonlyfont@\smc\fourteensmc \textonlyfont@\tt\fourteentt
%Erg"nzung des fetten Small-Capitals-Fonts:
%
 \ifsyntax@ \def\big##1{{\hbox{$\left##1\right.$}}}%
  \let\Big\big \let\bigg\big \let\Bigg\big
 \else
  \textfont\z@=\fourteenrm  \scriptfont\z@=\twelverm  \scriptscriptfont\z@=\tenrm
  \textfont\@ne=\fourteeni  \scriptfont\@ne=\twelvei  \scriptscriptfont\@ne=\teni
  \textfont\tw@=\fourteensy \scriptfont\tw@=\twelvesy \scriptscriptfont\tw@=\tensy
  \textfont\thr@@=\fourteenex \scriptfont\thr@@=\twelveex
        \scriptscriptfont\thr@@=\twelveex
  \textfont\itfam=\fourteenit \scriptfont\itfam=\twelveit
        \scriptscriptfont\itfam=\twelveit
  \textfont\bffam=\fourteenbf \scriptfont\bffam=\twelvebf
        \scriptscriptfont\bffam=\tenbf
  \setbox\strutbox\hbox{\vrule height12.2\p@ depth5\p@ width\z@}%
  \setbox\strutbox@\hbox{\lower.72\normallineskiplimit\vbox{%
        \kern-\normallineskiplimit\copy\strutbox}}%
 \setbox\z@\vbox{\hbox{$($}\kern\z@}\bigsize@=1.7\ht\z@
 \fi
 \normalbaselines\rm\ex@.2326ex\jot4.3\ex@\the\fourteenpoint@}

\font@\seventeenrm=cmr10 scaled\magstep3
\font@\seventeenit=cmti10 scaled\magstep3
\font@\seventeensl=cmsl10 scaled\magstep3
\font@\seventeensmc=cmcsc10 scaled\magstep3
\font@\seventeentt=cmtt10 scaled\magstep3
\font@\seventeenbf=cmbx10 scaled\magstep3
\font@\seventeeni=cmmi10 scaled\magstep3
\font@\seventeensy=cmsy10 scaled\magstep3
\font@\seventeenex=cmex10 scaled\magstep3
\font@\seventeenmsa=msam10 scaled\magstep3
\font@\seventeeneufm=eufm10 scaled\magstep3
\font@\seventeenmsb=msbm10 scaled\magstep3
\newtoks\seventeenpoint@
\def\seventeenpoint{\normalbaselineskip18\p@
 \abovedisplayskip21.6\p@ plus5.2\p@ minus15.4\p@
 \belowdisplayskip\abovedisplayskip
 \abovedisplayshortskip\z@ plus5.2\p@
 \belowdisplayshortskip12.1\p@ plus5.2\p@ minus7\p@
 \textonlyfont@\rm\seventeenrm \textonlyfont@\it\seventeenit
 \textonlyfont@\sl\seventeensl \textonlyfont@\bf\seventeenbf
 \textonlyfont@\smc\seventeensmc \textonlyfont@\tt\seventeentt
%Erg"nzung des fetten Small-Capitals-Fonts:
%
 \ifsyntax@ \def\big##1{{\hbox{$\left##1\right.$}}}%
  \let\Big\big \let\bigg\big \let\Bigg\big
 \else
  \textfont\z@=\seventeenrm  \scriptfont\z@=\fourteenrm  \scriptscriptfont\z@=\twelverm
  \textfont\@ne=\seventeeni  \scriptfont\@ne=\fourteeni  \scriptscriptfont\@ne=\twelvei
  \textfont\tw@=\seventeensy \scriptfont\tw@=\fourteensy \scriptscriptfont\tw@=\twelvesy
  \textfont\thr@@=\seventeenex \scriptfont\thr@@=\fourteenex
        \scriptscriptfont\thr@@=\fourteenex
  \textfont\itfam=\seventeenit \scriptfont\itfam=\fourteenit
        \scriptscriptfont\itfam=\fourteenit
  \textfont\bffam=\seventeenbf \scriptfont\bffam=\fourteenbf
        \scriptscriptfont\bffam=\twelvebf
  \setbox\strutbox\hbox{\vrule height14.6\p@ depth6\p@ width\z@}%
  \setbox\strutbox@\hbox{\lower.86\normallineskiplimit\vbox{%
        \kern-\normallineskiplimit\copy\strutbox}}%
 \setbox\z@\vbox{\hbox{$($}\kern\z@}\bigsize@=2\ht\z@
 \fi
 \normalbaselines\rm\ex@.2326ex\jot5.2\ex@\the\seventeenpoint@}

\catcode`\@=13
\catcode`\@=11
\font\tenln    = line10
\font\tenlnw   = linew10

\newskip\Einheit \Einheit=0.5cm
\newcount\xcoord \newcount\ycoord
\newdimen\xdim \newdimen\ydim \newdimen\PfadD@cke \newdimen\Pfadd@cke

%%%%%%%%%%%%%%%%%%%%%%%%%%%%%%%%%%%%%%%%%%%%%%%%%
%LaTeX counters, dimensions, variables for lines%
%%%%%%%%%%%%%%%%%%%%%%%%%%%%%%%%%%%%%%%%%%%%%%%%%
\newcount\@tempcnta
\newcount\@tempcntb

\newdimen\@tempdima
\newdimen\@tempdimb

\newdimen\@wholewidth
\newdimen\@halfwidth

\newcount\@xarg
\newcount\@yarg
\newcount\@yyarg
\newbox\@linechar
\newbox\@tempboxa
\newdimen\@linelen
\newdimen\@clnwd
\newdimen\@clnht

\newif\if@negarg

\def\@whilenoop#1{}
\def\@whiledim#1\do #2{\ifdim #1\relax#2\@iwhiledim{#1\relax#2}\fi}
\def\@iwhiledim#1{\ifdim #1\let\@nextwhile=\@iwhiledim
        \else\let\@nextwhile=\@whilenoop\fi\@nextwhile{#1}}

\def\@whileswnoop#1\fi{}
\def\@whilesw#1\fi#2{#1#2\@iwhilesw{#1#2}\fi\fi}
\def\@iwhilesw#1\fi{#1\let\@nextwhile=\@iwhilesw
         \else\let\@nextwhile=\@whileswnoop\fi\@nextwhile{#1}\fi}

\def\thinlines{\let\@linefnt\tenln \let\@circlefnt\tencirc
  \@wholewidth\fontdimen8\tenln \@halfwidth .5\@wholewidth}
\def\thicklines{\let\@linefnt\tenlnw \let\@circlefnt\tencircw
  \@wholewidth\fontdimen8\tenlnw \@halfwidth .5\@wholewidth}
\thinlines
%%%%%%%%%%%%%%%%%%%%%%%%%%%%%%%%%%%%%%%%%%%%%%%%%%%%%%%%%%%

\PfadD@cke1pt \Pfadd@cke0.5pt
\def\PfadDicke#1{\PfadD@cke#1 \divide\PfadD@cke by2 \Pfadd@cke\PfadD@cke \multiply\PfadD@cke by2}
\long\def\LOOP#1\REPEAT{\def\BODY{#1}\ITERATE}
\def\ITERATE{\BODY \let\next\ITERATE \else\let\next\relax\fi \next}
\let\REPEAT=\fi
\def\Punkt{\hbox{\raise-2pt\hbox to0pt{\hss$\ssize\bullet$\hss}}}
\def\DuennPunkt(#1,#2){\unskip
  \raise#2 \Einheit\hbox to0pt{\hskip#1 \Einheit
          \raise-2.5pt\hbox to0pt{\hss$\bullet$\hss}\hss}}
\def\NormalPunkt(#1,#2){\unskip
  \raise#2 \Einheit\hbox to0pt{\hskip#1 \Einheit
          \raise-3pt\hbox to0pt{\hss\twelvepoint$\bullet$\hss}\hss}}
\def\DickPunkt(#1,#2){\unskip
  \raise#2 \Einheit\hbox to0pt{\hskip#1 \Einheit
          \raise-4pt\hbox to0pt{\hss\fourteenpoint$\bullet$\hss}\hss}}
\def\Kreis(#1,#2){\unskip
  \raise#2 \Einheit\hbox to0pt{\hskip#1 \Einheit
          \raise-4pt\hbox to0pt{\hss\fourteenpoint$\circ$\hss}\hss}}

%%%%%%%%%%%%%%%%%%%%%
%LaTeX line macros%
%%%%%%%%%%%%%%%%%%%%%
\def\Line@(#1,#2)#3{\@xarg #1\relax \@yarg #2\relax
\@linelen=#3\Einheit
\ifnum\@xarg =0 \@vline
  \else \ifnum\@yarg =0 \@hline \else \@sline\fi
\fi}

\def\@sline{\ifnum\@xarg< 0 \@negargtrue \@xarg -\@xarg \@yyarg -\@yarg
  \else \@negargfalse \@yyarg \@yarg \fi
\ifnum \@yyarg >0 \@tempcnta\@yyarg \else \@tempcnta -\@yyarg \fi
\ifnum\@tempcnta>6 \@badlinearg\@tempcnta0 \fi
\ifnum\@xarg>6 \@badlinearg\@xarg 1 \fi
\setbox\@linechar\hbox{\@linefnt\@getlinechar(\@xarg,\@yyarg)}%
\ifnum \@yarg >0 \let\@upordown\raise \@clnht\z@
   \else\let\@upordown\lower \@clnht \ht\@linechar\fi
\@clnwd=\wd\@linechar
\if@negarg \hskip -\wd\@linechar \def\@tempa{\hskip -2\wd\@linechar}\else
     \let\@tempa\relax \fi
\@whiledim \@clnwd <\@linelen \do
  {\@upordown\@clnht\copy\@linechar
   \@tempa
   \advance\@clnht \ht\@linechar
   \advance\@clnwd \wd\@linechar}%
\advance\@clnht -\ht\@linechar
\advance\@clnwd -\wd\@linechar
\@tempdima\@linelen\advance\@tempdima -\@clnwd
\@tempdimb\@tempdima\advance\@tempdimb -\wd\@linechar
\if@negarg \hskip -\@tempdimb \else \hskip \@tempdimb \fi
\multiply\@tempdima \@m
\@tempcnta \@tempdima \@tempdima \wd\@linechar \divide\@tempcnta \@tempdima
\@tempdima \ht\@linechar \multiply\@tempdima \@tempcnta
\divide\@tempdima \@m
\advance\@clnht \@tempdima
\ifdim \@linelen <\wd\@linechar
   \hskip \wd\@linechar
  \else\@upordown\@clnht\copy\@linechar\fi}

\def\@hline{\ifnum \@xarg <0 \hskip -\@linelen \fi
\vrule height\Pfadd@cke width \@linelen depth\Pfadd@cke
\ifnum \@xarg <0 \hskip -\@linelen \fi}

\def\@getlinechar(#1,#2){\@tempcnta#1\relax\multiply\@tempcnta 8
\advance\@tempcnta -9 \ifnum #2>0 \advance\@tempcnta #2\relax\else
\advance\@tempcnta -#2\relax\advance\@tempcnta 64 \fi
\char\@tempcnta}

\def\Vektor(#1,#2)#3(#4,#5){\unskip\leavevmode
  \xcoord#4\relax \ycoord#5\relax
      \raise\ycoord \Einheit\hbox to0pt{\hskip\xcoord \Einheit
         \Vector@(#1,#2){#3}\hss}}

\def\Vector@(#1,#2)#3{\@xarg #1\relax \@yarg #2\relax
\@tempcnta \ifnum\@xarg<0 -\@xarg\else\@xarg\fi
\ifnum\@tempcnta<5\relax
\@linelen=#3\Einheit
\ifnum\@xarg =0 \@vvector
  \else \ifnum\@yarg =0 \@hvector \else \@svector\fi
\fi
\else\@badlinearg\fi}

\def\@hvector{\@hline\hbox to 0pt{\@linefnt
\ifnum \@xarg <0 \@getlarrow(1,0)\hss\else
    \hss\@getrarrow(1,0)\fi}}

\def\@vvector{\ifnum \@yarg <0 \@downvector \else \@upvector \fi}

\def\@svector{\@sline
\@tempcnta\@yarg \ifnum\@tempcnta <0 \@tempcnta=-\@tempcnta\fi
\ifnum\@tempcnta <5
  \hskip -\wd\@linechar
  \@upordown\@clnht \hbox{\@linefnt  \if@negarg
  \@getlarrow(\@xarg,\@yyarg) \else \@getrarrow(\@xarg,\@yyarg) \fi}%
\else\@badlinearg\fi}

\def\@upline{\hbox to \z@{\hskip -.5\Pfadd@cke \vrule width \Pfadd@cke
   height \@linelen depth \z@\hss}}

\def\@downline{\hbox to \z@{\hskip -.5\Pfadd@cke \vrule width \Pfadd@cke
   height \z@ depth \@linelen \hss}}

\def\@upvector{\@upline\setbox\@tempboxa\hbox{\@linefnt\char'66}\raise
     \@linelen \hbox to\z@{\lower \ht\@tempboxa\box\@tempboxa\hss}}

\def\@downvector{\@downline\lower \@linelen
      \hbox to \z@{\@linefnt\char'77\hss}}

\def\@getlarrow(#1,#2){\ifnum #2 =\z@ \@tempcnta='33\else
\@tempcnta=#1\relax\multiply\@tempcnta \sixt@@n \advance\@tempcnta
-9 \@tempcntb=#2\relax\multiply\@tempcntb \tw@
\ifnum \@tempcntb >0 \advance\@tempcnta \@tempcntb\relax
\else\advance\@tempcnta -\@tempcntb\advance\@tempcnta 64
\fi\fi\char\@tempcnta}

\def\@getrarrow(#1,#2){\@tempcntb=#2\relax
\ifnum\@tempcntb < 0 \@tempcntb=-\@tempcntb\relax\fi
\ifcase \@tempcntb\relax \@tempcnta='55 \or
\ifnum #1<3 \@tempcnta=#1\relax\multiply\@tempcnta
24 \advance\@tempcnta -6 \else \ifnum #1=3 \@tempcnta=49
\else\@tempcnta=58 \fi\fi\or
\ifnum #1<3 \@tempcnta=#1\relax\multiply\@tempcnta
24 \advance\@tempcnta -3 \else \@tempcnta=51\fi\or
\@tempcnta=#1\relax\multiply\@tempcnta
\sixt@@n \advance\@tempcnta -\tw@ \else
\@tempcnta=#1\relax\multiply\@tempcnta
\sixt@@n \advance\@tempcnta 7 \fi\ifnum #2<0 \advance\@tempcnta 64 \fi
\char\@tempcnta}
%%%%%%%%%%%%%%%%%%%%%%%%%%%%%%%%%%%%%%%%%%%%%%%%%%%%%%%%%%%%%

\def\Diagonale(#1,#2)#3{\unskip\leavevmode
  \xcoord#1\relax \ycoord#2\relax
      \raise\ycoord \Einheit\hbox to0pt{\hskip\xcoord \Einheit
         \Line@(1,1){#3}\hss}}
\def\AntiDiagonale(#1,#2)#3{\unskip\leavevmode
  \xcoord#1\relax \ycoord#2\relax %\advance\xcoord by -0.05\relax
      \raise\ycoord \Einheit\hbox to0pt{\hskip\xcoord \Einheit
         \Line@(1,-1){#3}\hss}}
\def\Pfad(#1,#2),#3\endPfad{\unskip\leavevmode
  \xcoord#1 \ycoord#2 \thicklines\ZeichnePfad#3\endPfad\thinlines}
\def\ZeichnePfad#1{\ifx#1\endPfad\let\next\relax
  \else\let\next\ZeichnePfad
    \ifnum#1=1
      \raise\ycoord \Einheit\hbox to0pt{\hskip\xcoord \Einheit
         \vrule height\Pfadd@cke width1 \Einheit depth\Pfadd@cke\hss}%
      \advance\xcoord by 1
    \else\ifnum#1=2
      \raise\ycoord \Einheit\hbox to0pt{\hskip\xcoord \Einheit
        \hbox{\hskip-\PfadD@cke\vrule height1 \Einheit width\PfadD@cke depth0pt}\hss}%
      \advance\ycoord by 1
    \else\ifnum#1=3
      \raise\ycoord \Einheit\hbox to0pt{\hskip\xcoord \Einheit
         \Line@(1,1){1}\hss}
      \advance\xcoord by 1
      \advance\ycoord by 1
    \else\ifnum#1=4
      \raise\ycoord \Einheit\hbox to0pt{\hskip\xcoord \Einheit
         \Line@(1,-1){1}\hss}
      \advance\xcoord by 1
      \advance\ycoord by -1
    \fi\fi\fi\fi
  \fi\next}
\def\hSSchritt{\leavevmode\raise-.4pt\hbox to0pt{\hss.\hss}\hskip.2\Einheit
  \raise-.4pt\hbox to0pt{\hss.\hss}\hskip.2\Einheit
  \raise-.4pt\hbox to0pt{\hss.\hss}\hskip.2\Einheit
  \raise-.4pt\hbox to0pt{\hss.\hss}\hskip.2\Einheit
  \raise-.4pt\hbox to0pt{\hss.\hss}\hskip.2\Einheit}
\def\vSSchritt{\vbox{\baselineskip.2\Einheit\lineskiplimit0pt
\hbox{.}\hbox{.}\hbox{.}\hbox{.}\hbox{.}}}
\def\DSSchritt{\leavevmode\raise-.4pt\hbox to0pt{%
  \hbox to0pt{\hss.\hss}\hskip.2\Einheit
  \raise.2\Einheit\hbox to0pt{\hss.\hss}\hskip.2\Einheit
  \raise.4\Einheit\hbox to0pt{\hss.\hss}\hskip.2\Einheit
  \raise.6\Einheit\hbox to0pt{\hss.\hss}\hskip.2\Einheit
  \raise.8\Einheit\hbox to0pt{\hss.\hss}\hss}}
\def\dSSchritt{\leavevmode\raise-.4pt\hbox to0pt{%
  \hbox to0pt{\hss.\hss}\hskip.2\Einheit
  \raise-.2\Einheit\hbox to0pt{\hss.\hss}\hskip.2\Einheit
  \raise-.4\Einheit\hbox to0pt{\hss.\hss}\hskip.2\Einheit
  \raise-.6\Einheit\hbox to0pt{\hss.\hss}\hskip.2\Einheit
  \raise-.8\Einheit\hbox to0pt{\hss.\hss}\hss}}
\def\SPfad(#1,#2),#3\endSPfad{\unskip\leavevmode
  \xcoord#1 \ycoord#2 \ZeichneSPfad#3\endSPfad}
\def\ZeichneSPfad#1{\ifx#1\endSPfad\let\next\relax
  \else\let\next\ZeichneSPfad
    \ifnum#1=1
      \raise\ycoord \Einheit\hbox to0pt{\hskip\xcoord \Einheit
         \hSSchritt\hss}%
      \advance\xcoord by 1
    \else\ifnum#1=2
      \raise\ycoord \Einheit\hbox to0pt{\hskip\xcoord \Einheit
        \hbox{\hskip-2pt \vSSchritt}\hss}%
      \advance\ycoord by 1
    \else\ifnum#1=3
      \raise\ycoord \Einheit\hbox to0pt{\hskip\xcoord \Einheit
         \DSSchritt\hss}
      \advance\xcoord by 1
      \advance\ycoord by 1
    \else\ifnum#1=4
      \raise\ycoord \Einheit\hbox to0pt{\hskip\xcoord \Einheit
         \dSSchritt\hss}
      \advance\xcoord by 1
      \advance\ycoord by -1
    \fi\fi\fi\fi
  \fi\next}
\def\Koordinatenachsen(#1,#2){\unskip
 \hbox to0pt{\hskip-.5pt\vrule height#2 \Einheit width.5pt depth1 \Einheit}%
 \hbox to0pt{\hskip-1 \Einheit \xcoord#1 \advance\xcoord by1
    \vrule height0.25pt width\xcoord \Einheit depth0.25pt\hss}}
\def\Koordinatenachsen(#1,#2)(#3,#4){\unskip
 \hbox to0pt{\hskip-.5pt \ycoord-#4 \advance\ycoord by1
    \vrule height#2 \Einheit width.5pt depth\ycoord \Einheit}%
 \hbox to0pt{\hskip-1 \Einheit \hskip#3\Einheit
    \xcoord#1 \advance\xcoord by1 \advance\xcoord by-#3
    \vrule height0.25pt width\xcoord \Einheit depth0.25pt\hss}}
\def\Gitter(#1,#2){\unskip \xcoord0 \ycoord0 \leavevmode
  \LOOP\ifnum\ycoord<#2
    \loop\ifnum\xcoord<#1
      \raise\ycoord \Einheit\hbox to0pt{\hskip\xcoord \Einheit\Punkt\hss}%
      \advance\xcoord by1
    \repeat
    \xcoord0
    \advance\ycoord by1
  \REPEAT}
\def\Gitter(#1,#2)(#3,#4){\unskip \xcoord#3 \ycoord#4 \leavevmode
  \LOOP\ifnum\ycoord<#2
    \loop\ifnum\xcoord<#1
      \raise\ycoord \Einheit\hbox to0pt{\hskip\xcoord \Einheit\Punkt\hss}%
      \advance\xcoord by1
    \repeat
    \xcoord#3
    \advance\ycoord by1
  \REPEAT}
\def\Label#1#2(#3,#4){\unskip \xdim#3 \Einheit \ydim#4 \Einheit
  \def\lo{\advance\xdim by-.5 \Einheit \advance\ydim by.5 \Einheit}%
  \def\llo{\advance\xdim by-.25cm \advance\ydim by.5 \Einheit}%
  \def\loo{\advance\xdim by-.5 \Einheit \advance\ydim by.25cm}%
  \def\o{\advance\ydim by.25cm}%
  \def\ro{\advance\xdim by.5 \Einheit \advance\ydim by.5 \Einheit}%
  \def\rro{\advance\xdim by.25cm \advance\ydim by.5 \Einheit}%
  \def\roo{\advance\xdim by.5 \Einheit \advance\ydim by.25cm}%
  \def\l{\advance\xdim by-.30cm}%
  \def\r{\advance\xdim by.30cm}%
  \def\lu{\advance\xdim by-.5 \Einheit \advance\ydim by-.6 \Einheit}%
  \def\llu{\advance\xdim by-.25cm \advance\ydim by-.6 \Einheit}%
  \def\luu{\advance\xdim by-.5 \Einheit \advance\ydim by-.30cm}%
  \def\u{\advance\ydim by-.30cm}%
  \def\ru{\advance\xdim by.5 \Einheit \advance\ydim by-.6 \Einheit}%
  \def\rru{\advance\xdim by.25cm \advance\ydim by-.6 \Einheit}%
  \def\ruu{\advance\xdim by.5 \Einheit \advance\ydim by-.30cm}%
  #1\raise\ydim\hbox to0pt{\hskip\xdim
     \vbox to0pt{\vss\hbox to0pt{\hss$#2$\hss}\vss}\hss}%
}
\catcode`\@=13

\catcode`\@=11
\font\tenln    = line10
\font\tenlnw   = linew10

\newskip\Einheit \Einheit=0.5cm
\newcount\xcoord \newcount\ycoord
\newdimen\xdim \newdimen\ydim \newdimen\PfadD@cke \newdimen\Pfadd@cke

%%%%%%%%%%%%%%%%%%%%%%%%%%%%%%%%%%%%%%%%%%%%%%%%%
%LaTeX counters, dimensions, variables for lines%
%%%%%%%%%%%%%%%%%%%%%%%%%%%%%%%%%%%%%%%%%%%%%%%%%
\newcount\@tempcnta
\newcount\@tempcntb

\newdimen\@tempdima
\newdimen\@tempdimb

\newdimen\@wholewidth
\newdimen\@halfwidth

\newcount\@xarg
\newcount\@yarg
\newcount\@yyarg
\newbox\@linechar
\newbox\@tempboxa
\newdimen\@linelen
\newdimen\@clnwd
\newdimen\@clnht

\newif\if@negarg

\def\@whilenoop#1{}
\def\@whiledim#1\do #2{\ifdim #1\relax#2\@iwhiledim{#1\relax#2}\fi}
\def\@iwhiledim#1{\ifdim #1\let\@nextwhile=\@iwhiledim
        \else\let\@nextwhile=\@whilenoop\fi\@nextwhile{#1}}

\def\@whileswnoop#1\fi{}
\def\@whilesw#1\fi#2{#1#2\@iwhilesw{#1#2}\fi\fi}
\def\@iwhilesw#1\fi{#1\let\@nextwhile=\@iwhilesw
         \else\let\@nextwhile=\@whileswnoop\fi\@nextwhile{#1}\fi}

\def\thinlines{\let\@linefnt\tenln \let\@circlefnt\tencirc
  \@wholewidth\fontdimen8\tenln \@halfwidth .5\@wholewidth}
\def\thicklines{\let\@linefnt\tenlnw \let\@circlefnt\tencircw
  \@wholewidth\fontdimen8\tenlnw \@halfwidth .5\@wholewidth}
\thinlines
%%%%%%%%%%%%%%%%%%%%%%%%%%%%%%%%%%%%%%%%%%%%%%%%%%%%%%%%%%%

\PfadD@cke1pt \Pfadd@cke0.5pt
\def\PfadDicke#1{\PfadD@cke#1 \divide\PfadD@cke by2 \Pfadd@cke\PfadD@cke \multiply\PfadD@cke by2}
\long\def\LOOP#1\REPEAT{\def\BODY{#1}\ITERATE}
\def\ITERATE{\BODY \let\next\ITERATE \else\let\next\relax\fi \next}
\let\REPEAT=\fi
\def\Punkt{\hbox{\raise-2pt\hbox to0pt{\hss$\ssize\bullet$\hss}}}
\def\DuennPunkt(#1,#2){\unskip
  \raise#2 \Einheit\hbox to0pt{\hskip#1 \Einheit
          \raise-2.5pt\hbox to0pt{\hss$\bullet$\hss}\hss}}
\def\NormalPunkt(#1,#2){\unskip
  \raise#2 \Einheit\hbox to0pt{\hskip#1 \Einheit
          \raise-3pt\hbox to0pt{\hss\twelvepoint$\bullet$\hss}\hss}}
\def\DickPunkt(#1,#2){\unskip
  \raise#2 \Einheit\hbox to0pt{\hskip#1 \Einheit
          \raise-4pt\hbox to0pt{\hss\fourteenpoint$\bullet$\hss}\hss}}
\def\Kreis(#1,#2){\unskip
  \raise#2 \Einheit\hbox to0pt{\hskip#1 \Einheit
          \raise-4pt\hbox to0pt{\hss\fourteenpoint$\circ$\hss}\hss}}

%%%%%%%%%%%%%%%%%%%%%
%LaTeX line macros%
%%%%%%%%%%%%%%%%%%%%%
\def\Line@(#1,#2)#3{\@xarg #1\relax \@yarg #2\relax
\@linelen=#3\Einheit
\ifnum\@xarg =0 \@vline
  \else \ifnum\@yarg =0 \@hline \else \@sline\fi
\fi}

\def\@sline{\ifnum\@xarg< 0 \@negargtrue \@xarg -\@xarg \@yyarg -\@yarg
  \else \@negargfalse \@yyarg \@yarg \fi
\ifnum \@yyarg >0 \@tempcnta\@yyarg \else \@tempcnta -\@yyarg \fi
\ifnum\@tempcnta>6 \@badlinearg\@tempcnta0 \fi
\ifnum\@xarg>6 \@badlinearg\@xarg 1 \fi
\setbox\@linechar\hbox{\@linefnt\@getlinechar(\@xarg,\@yyarg)}%
\ifnum \@yarg >0 \let\@upordown\raise \@clnht\z@
   \else\let\@upordown\lower \@clnht \ht\@linechar\fi
\@clnwd=\wd\@linechar
\if@negarg \hskip -\wd\@linechar \def\@tempa{\hskip -2\wd\@linechar}\else
     \let\@tempa\relax \fi
\@whiledim \@clnwd <\@linelen \do
  {\@upordown\@clnht\copy\@linechar
   \@tempa
   \advance\@clnht \ht\@linechar
   \advance\@clnwd \wd\@linechar}%
\advance\@clnht -\ht\@linechar
\advance\@clnwd -\wd\@linechar
\@tempdima\@linelen\advance\@tempdima -\@clnwd
\@tempdimb\@tempdima\advance\@tempdimb -\wd\@linechar
\if@negarg \hskip -\@tempdimb \else \hskip \@tempdimb \fi
\multiply\@tempdima \@m
\@tempcnta \@tempdima \@tempdima \wd\@linechar \divide\@tempcnta \@tempdima
\@tempdima \ht\@linechar \multiply\@tempdima \@tempcnta
\divide\@tempdima \@m
\advance\@clnht \@tempdima
\ifdim \@linelen <\wd\@linechar
   \hskip \wd\@linechar
  \else\@upordown\@clnht\copy\@linechar\fi}

\def\@hline{\ifnum \@xarg <0 \hskip -\@linelen \fi
\vrule height\Pfadd@cke width \@linelen depth\Pfadd@cke
\ifnum \@xarg <0 \hskip -\@linelen \fi}

\def\@getlinechar(#1,#2){\@tempcnta#1\relax\multiply\@tempcnta 8
\advance\@tempcnta -9 \ifnum #2>0 \advance\@tempcnta #2\relax\else
\advance\@tempcnta -#2\relax\advance\@tempcnta 64 \fi
\char\@tempcnta}

\def\Vektor(#1,#2)#3(#4,#5){\unskip\leavevmode
  \xcoord#4\relax \ycoord#5\relax
      \raise\ycoord \Einheit\hbox to0pt{\hskip\xcoord \Einheit
         \Vector@(#1,#2){#3}\hss}}

\def\Vector@(#1,#2)#3{\@xarg #1\relax \@yarg #2\relax
\@tempcnta \ifnum\@xarg<0 -\@xarg\else\@xarg\fi
\ifnum\@tempcnta<5\relax
\@linelen=#3\Einheit
\ifnum\@xarg =0 \@vvector
  \else \ifnum\@yarg =0 \@hvector \else \@svector\fi
\fi
\else\@badlinearg\fi}

\def\@hvector{\@hline\hbox to 0pt{\@linefnt
\ifnum \@xarg <0 \@getlarrow(1,0)\hss\else
    \hss\@getrarrow(1,0)\fi}}

\def\@vvector{\ifnum \@yarg <0 \@downvector \else \@upvector \fi}

\def\@svector{\@sline
\@tempcnta\@yarg \ifnum\@tempcnta <0 \@tempcnta=-\@tempcnta\fi
\ifnum\@tempcnta <5
  \hskip -\wd\@linechar
  \@upordown\@clnht \hbox{\@linefnt  \if@negarg
  \@getlarrow(\@xarg,\@yyarg) \else \@getrarrow(\@xarg,\@yyarg) \fi}%
\else\@badlinearg\fi}

\def\@upline{\hbox to \z@{\hskip -.5\Pfadd@cke \vrule width \Pfadd@cke
   height \@linelen depth \z@\hss}}

\def\@downline{\hbox to \z@{\hskip -.5\Pfadd@cke \vrule width \Pfadd@cke
   height \z@ depth \@linelen \hss}}

\def\@upvector{\@upline\setbox\@tempboxa\hbox{\@linefnt\char'66}\raise
     \@linelen \hbox to\z@{\lower \ht\@tempboxa\box\@tempboxa\hss}}

\def\@downvector{\@downline\lower \@linelen
      \hbox to \z@{\@linefnt\char'77\hss}}

\def\@getlarrow(#1,#2){\ifnum #2 =\z@ \@tempcnta='33\else
\@tempcnta=#1\relax\multiply\@tempcnta \sixt@@n \advance\@tempcnta
-9 \@tempcntb=#2\relax\multiply\@tempcntb \tw@
\ifnum \@tempcntb >0 \advance\@tempcnta \@tempcntb\relax
\else\advance\@tempcnta -\@tempcntb\advance\@tempcnta 64
\fi\fi\char\@tempcnta}

\def\@getrarrow(#1,#2){\@tempcntb=#2\relax
\ifnum\@tempcntb < 0 \@tempcntb=-\@tempcntb\relax\fi
\ifcase \@tempcntb\relax \@tempcnta='55 \or
\ifnum #1<3 \@tempcnta=#1\relax\multiply\@tempcnta
24 \advance\@tempcnta -6 \else \ifnum #1=3 \@tempcnta=49
\else\@tempcnta=58 \fi\fi\or
\ifnum #1<3 \@tempcnta=#1\relax\multiply\@tempcnta
24 \advance\@tempcnta -3 \else \@tempcnta=51\fi\or
\@tempcnta=#1\relax\multiply\@tempcnta
\sixt@@n \advance\@tempcnta -\tw@ \else
\@tempcnta=#1\relax\multiply\@tempcnta
\sixt@@n \advance\@tempcnta 7 \fi\ifnum #2<0 \advance\@tempcnta 64 \fi
\char\@tempcnta}
%%%%%%%%%%%%%%%%%%%%%%%%%%%%%%%%%%%%%%%%%%%%%%%%%%%%%%%%%%%%%

\def\Diagonale(#1,#2)#3{\unskip\leavevmode
  \xcoord#1\relax \ycoord#2\relax
      \raise\ycoord \Einheit\hbox to0pt{\hskip\xcoord \Einheit
         \Line@(1,1){#3}\hss}}
\def\AntiDiagonale(#1,#2)#3{\unskip\leavevmode
  \xcoord#1\relax \ycoord#2\relax %\advance\xcoord by -0.05\relax
      \raise\ycoord \Einheit\hbox to0pt{\hskip\xcoord \Einheit
         \Line@(1,-1){#3}\hss}}
\def\Pfad(#1,#2),#3\endPfad{\unskip\leavevmode
  \xcoord#1 \ycoord#2 \thicklines\ZeichnePfad#3\endPfad\thinlines}
\def\ZeichnePfad#1{\ifx#1\endPfad\let\next\relax
  \else\let\next\ZeichnePfad
    \ifnum#1=1
      \raise\ycoord \Einheit\hbox to0pt{\hskip\xcoord \Einheit
         \vrule height\Pfadd@cke width1 \Einheit depth\Pfadd@cke\hss}%
      \advance\xcoord by 1
    \else\ifnum#1=2
      \raise\ycoord \Einheit\hbox to0pt{\hskip\xcoord \Einheit
        \hbox{\hskip-\PfadD@cke\vrule height1 \Einheit width\PfadD@cke depth0pt}\hss}%
      \advance\ycoord by 1
    \else\ifnum#1=3
      \raise\ycoord \Einheit\hbox to0pt{\hskip\xcoord \Einheit
         \Line@(1,1){1}\hss}
      \advance\xcoord by 1
      \advance\ycoord by 1
    \else\ifnum#1=4
      \raise\ycoord \Einheit\hbox to0pt{\hskip\xcoord \Einheit
         \Line@(1,-1){1}\hss}
      \advance\xcoord by 1
      \advance\ycoord by -1
    \else\ifnum#1=5
      \advance\xcoord by -1
      \raise\ycoord \Einheit\hbox to0pt{\hskip\xcoord \Einheit
         \vrule height\Pfadd@cke width1 \Einheit depth\Pfadd@cke\hss}%
    \else\ifnum#1=6
      \advance\ycoord by -1
      \raise\ycoord \Einheit\hbox to0pt{\hskip\xcoord \Einheit
        \hbox{\hskip-\PfadD@cke\vrule height1 \Einheit width\PfadD@cke depth0pt}\hss}%
    \else\ifnum#1=7
      \advance\xcoord by -1
      \advance\ycoord by -1
      \raise\ycoord \Einheit\hbox to0pt{\hskip\xcoord \Einheit
         \Line@(1,1){1}\hss}
    \else\ifnum#1=8
      \advance\xcoord by -1
      \advance\ycoord by +1
      \raise\ycoord \Einheit\hbox to0pt{\hskip\xcoord \Einheit
         \Line@(1,-1){1}\hss}
    \fi\fi\fi\fi
    \fi\fi\fi\fi
  \fi\next}
\def\hSSchritt{\leavevmode\raise-.4pt\hbox to0pt{\hss.\hss}\hskip.2\Einheit
  \raise-.4pt\hbox to0pt{\hss.\hss}\hskip.2\Einheit
  \raise-.4pt\hbox to0pt{\hss.\hss}\hskip.2\Einheit
  \raise-.4pt\hbox to0pt{\hss.\hss}\hskip.2\Einheit
  \raise-.4pt\hbox to0pt{\hss.\hss}\hskip.2\Einheit}
\def\vSSchritt{\vbox{\baselineskip.2\Einheit\lineskiplimit0pt
\hbox{.}\hbox{.}\hbox{.}\hbox{.}\hbox{.}}}
\def\DSSchritt{\leavevmode\raise-.4pt\hbox to0pt{%
  \hbox to0pt{\hss.\hss}\hskip.2\Einheit
  \raise.2\Einheit\hbox to0pt{\hss.\hss}\hskip.2\Einheit
  \raise.4\Einheit\hbox to0pt{\hss.\hss}\hskip.2\Einheit
  \raise.6\Einheit\hbox to0pt{\hss.\hss}\hskip.2\Einheit
  \raise.8\Einheit\hbox to0pt{\hss.\hss}\hss}}
\def\dSSchritt{\leavevmode\raise-.4pt\hbox to0pt{%
  \hbox to0pt{\hss.\hss}\hskip.2\Einheit
  \raise-.2\Einheit\hbox to0pt{\hss.\hss}\hskip.2\Einheit
  \raise-.4\Einheit\hbox to0pt{\hss.\hss}\hskip.2\Einheit
  \raise-.6\Einheit\hbox to0pt{\hss.\hss}\hskip.2\Einheit
  \raise-.8\Einheit\hbox to0pt{\hss.\hss}\hss}}
\def\SPfad(#1,#2),#3\endSPfad{\unskip\leavevmode
  \xcoord#1 \ycoord#2 \ZeichneSPfad#3\endSPfad}
\def\ZeichneSPfad#1{\ifx#1\endSPfad\let\next\relax
  \else\let\next\ZeichneSPfad
    \ifnum#1=1
      \raise\ycoord \Einheit\hbox to0pt{\hskip\xcoord \Einheit
         \hSSchritt\hss}%
      \advance\xcoord by 1
    \else\ifnum#1=2
      \raise\ycoord \Einheit\hbox to0pt{\hskip\xcoord \Einheit
        \hbox{\hskip-2pt \vSSchritt}\hss}%
      \advance\ycoord by 1
    \else\ifnum#1=3
      \raise\ycoord \Einheit\hbox to0pt{\hskip\xcoord \Einheit
         \DSSchritt\hss}
      \advance\xcoord by 1
      \advance\ycoord by 1
    \else\ifnum#1=4
      \raise\ycoord \Einheit\hbox to0pt{\hskip\xcoord \Einheit
         \dSSchritt\hss}
      \advance\xcoord by 1
      \advance\ycoord by -1
    \else\ifnum#1=5
      \advance\xcoord by -1
      \raise\ycoord \Einheit\hbox to0pt{\hskip\xcoord \Einheit
         \hSSchritt\hss}%
    \else\ifnum#1=6
      \advance\ycoord by -1
      \raise\ycoord \Einheit\hbox to0pt{\hskip\xcoord \Einheit
        \hbox{\hskip-2pt \vSSchritt}\hss}%
    \else\ifnum#1=7
      \advance\xcoord by -1
      \advance\ycoord by -1
      \raise\ycoord \Einheit\hbox to0pt{\hskip\xcoord \Einheit
         \DSSchritt\hss}
    \else\ifnum#1=8
      \advance\xcoord by -1
      \advance\ycoord by 1
      \raise\ycoord \Einheit\hbox to0pt{\hskip\xcoord \Einheit
         \dSSchritt\hss}
    \fi\fi\fi\fi
    \fi\fi\fi\fi
  \fi\next}
\def\Koordinatenachsen(#1,#2){\unskip
 \hbox to0pt{\hskip-.5pt\vrule height#2 \Einheit width.5pt depth1 \Einheit}%
 \hbox to0pt{\hskip-1 \Einheit \xcoord#1 \advance\xcoord by1
    \vrule height0.25pt width\xcoord \Einheit depth0.25pt\hss}}
\def\Koordinatenachsen(#1,#2)(#3,#4){\unskip
 \hbox to0pt{\hskip-.5pt \ycoord-#4 \advance\ycoord by1
    \vrule height#2 \Einheit width.5pt depth\ycoord \Einheit}%
 \hbox to0pt{\hskip-1 \Einheit \hskip#3\Einheit
    \xcoord#1 \advance\xcoord by1 \advance\xcoord by-#3
    \vrule height0.25pt width\xcoord \Einheit depth0.25pt\hss}}
\def\Gitter(#1,#2){\unskip \xcoord0 \ycoord0 \leavevmode
  \LOOP\ifnum\ycoord<#2
    \loop\ifnum\xcoord<#1
      \raise\ycoord \Einheit\hbox to0pt{\hskip\xcoord \Einheit\Punkt\hss}%
      \advance\xcoord by1
    \repeat
    \xcoord0
    \advance\ycoord by1
  \REPEAT}
\def\Gitter(#1,#2)(#3,#4){\unskip \xcoord#3 \ycoord#4 \leavevmode
  \LOOP\ifnum\ycoord<#2
    \loop\ifnum\xcoord<#1
      \raise\ycoord \Einheit\hbox to0pt{\hskip\xcoord \Einheit\Punkt\hss}%
      \advance\xcoord by1
    \repeat
    \xcoord#3
    \advance\ycoord by1
  \REPEAT}
\def\Label#1#2(#3,#4){\unskip \xdim#3 \Einheit \ydim#4 \Einheit
  \def\lo{\advance\xdim by-.5 \Einheit \advance\ydim by.5 \Einheit}%
  \def\llo{\advance\xdim by-.25cm \advance\ydim by.5 \Einheit}%
  \def\loo{\advance\xdim by-.5 \Einheit \advance\ydim by.25cm}%
  \def\o{\advance\ydim by.25cm}%
  \def\ro{\advance\xdim by.5 \Einheit \advance\ydim by.5 \Einheit}%
  \def\rro{\advance\xdim by.25cm \advance\ydim by.5 \Einheit}%
  \def\roo{\advance\xdim by.5 \Einheit \advance\ydim by.25cm}%
  \def\l{\advance\xdim by-.30cm}%
  \def\r{\advance\xdim by.30cm}%
  \def\lu{\advance\xdim by-.5 \Einheit \advance\ydim by-.6 \Einheit}%
  \def\llu{\advance\xdim by-.25cm \advance\ydim by-.6 \Einheit}%
  \def\luu{\advance\xdim by-.5 \Einheit \advance\ydim by-.30cm}%
  \def\u{\advance\ydim by-.30cm}%
  \def\ru{\advance\xdim by.5 \Einheit \advance\ydim by-.6 \Einheit}%
  \def\rru{\advance\xdim by.25cm \advance\ydim by-.6 \Einheit}%
  \def\ruu{\advance\xdim by.5 \Einheit \advance\ydim by-.30cm}%
  #1\raise\ydim\hbox to0pt{\hskip\xdim
     \vbox to0pt{\vss\hbox to0pt{\hss$#2$\hss}\vss}\hss}%
}
\catcode`\@=13

\TagsOnRight

\def\RoZeAA{13}
\def\ModaAA{12}
\def\LindAA{11}
\def\KulkAC{10}
\def\KrLaAA{9}
\def\KrRuAA{8}
\def\KrPrAA{7}
\def\KratBQ{6}
\def\KratBE{5}
\def\GeViAB{4}
\def\FultAC{3}
\def\CoHeAA{2}
\def\BiLaAA{1}

\def\mult{\operatorname{mult}}
\def\EN{\operatorname{EN}}
\def\GF{\operatorname{GF}}
\def\Gr{\operatorname{Gr}}
\def\3{\ss}

\topmatter
\title On multiplicities of points on Schubert varieties in
Gra\3mannians II
\endtitle
\author C.~Krattenthaler\footnote"$^\dagger$"{Partially supported
by EC's IHRP Programme, grant RTN2-2001-00059, and by the
Austrian\linebreak
\hbox{Science Foundation FWF, grant P13190-MAT.}\hss}
\endauthor
\affil
Institut f\"ur Mathematik der Universit\"at Wien,\\
Strudlhofgasse 4, A-1090 Wien, Austria.\\
e-mail: KRATT\@Ap.Univie.Ac.At\\
WWW: \tt http://www.mat.univie.ac.at/\~{}kratt
\endaffil
\address Institut f\"ur Mathematik der Universit\"at Wien,
Strudlhofgasse 4, A-1090 Wien, Austria.
\endaddress
%\email KRATT\@Ap.Univie.Ac.At\\
%WWW: \tt http://www.mat.univie.ac.at/People/kratt\endemail
%\dedicatory \enddedicatory
%\date \enddate
%\thanks \endthanks
\abstract
We prove a conjecture by Kreiman and Lakshmibai on a combinatorial
description of multiplicities of points on Schubert
varieties in Gra\3mannians in terms of certain sets of reflections in
the corresponding Weyl group. The proof is accomplished by setting
up a bijection between these sets of reflections and the author's previous
combinatorial interpretation of these multiplicities in terms of
nonintersecting lattice paths (S\'eminaire Lotharingien Combin\.
{\bf 45} (2001), Article~B45c).
\endabstract
\endtopmatter
\document

\leftheadtext{C. Krattenthaler}

\subhead 1. Introduction\endsubhead
The {\it multiplicity} of a point on an algebraic variety is an important
invariant that ``measures" singularity of the point.
It was an important recent advance in Schubert calculus when
Rosenthal and Zelevinsky \cite{\RoZeAA} gave a
determinantal formula for the multiplicity of a point on
a Schubert variety in a Gra\3mannian. It paved the way to a
combinatorial understanding of this multiplicity. More precisely, it
was shown in \cite{\KratBQ} that it counts certain families of
nonintersecting lattice paths (and also certain tableaux). An
alternative, conjectural
combinatorial interpretation was proposed by Kreiman and
Lakshmibai in \cite{\KrLaAA, Conjecture~2}, in terms of certain sets of
reflections. The purpose of this paper is to prove that this latter
combinatorial interpretation is indeed valid.

The reason for the proposition of this alternative combinatorial
interpretation of the multiplicity of a point on
a Schubert variety in a Gra\3mannian in terms of sets of reflections
is that it appears that these sets of reflections also allow the
computation of the Hilbert series of the tangent cone at this
point (see \cite{\KrLaAA, Conjecture~1}).
While we are not able to prove this more general conjecture, we provide an
equivalent form of the conjecture in which the Hilbert series is
essentially given in terms of a generating function for certain
families of nonintersecting lattice paths which are counted with
respect to turns. This equivalent form of the conjecture has the
advantage over the original form that it reduces the computation of
the Hilbert series to a finite problem. Moreover, it is analogous to
similar formulas for the Hilbert series associated to related
determinantal varieties (see \cite{\CoHeAA, Eq.~(1)},
\cite{\KrPrAA, p.~1021, line~11}
or \cite{\KrRuAA, Theorem~1}).

Our paper is organised as follows.
In the next section we fix notation and formulate the multiplicity
conjecture by Kreiman and Lakshmibai. There we also recall the
author's combinatorial interpretation of the multiplicity in terms of
nonintersecting lattice paths. Section~3 contains the proof of the
conjecture, which is accomplished
by setting up a bijection between these nonintersecting
lattice paths and the sets of reflections of Kreiman and Lakshmibai.
Finally, in Theorem~2 in Section~4, we show that the results from Section~3
allow in fact the above mentioned reformulation, in terms of
nonintersecting lattice paths, of
Conjecture~1 in \cite{\KrLaAA} on the Hilbert series of the tangent
cone at a point on a Schubert variety in the Gra\3mannian.
It is an open problem to find a compact formula for the generating
function of nonintersecting lattice paths that appears in this
formulation (see Remark~(2) after Theorem~2).

\subhead 2. Combinatorial interpretations of multiplicities of points
on Schubert varieties in Gra\3mannians\endsubhead
We recall some basic notions from the Schubert calculus
in the Gra\3mannian, and fix the notation that we are going to use.
We refer the reader to \cite{\BiLaAA, Sec.~3.1} and \cite{\FultAC, Sec.~9.4}
for in-depth introductions into the subject.

\smallskip
Let $d$ and $n$ be positive integers with $0\le d\le n$. The {\it
Gra\3mannian} $\Gr_d(V)$ is the variety of all $d$-dimensional
subspaces in an $n$-dimensional vector space $V$ (over some
algebraically closed field of arbitrary characteristic).
Schubert varieties in the Gra\3mannian $\Gr_d(V)$ are indexed by
elements in
$S_n/(S_d\times S_{n-d})$, where $S_m$ denotes the symmetric group of
order $m$. Any coset $C$ in $S_n/(S_d\times S_{n-d})$ has a {\it minimal
representative}, which is the unique permutation $w=i_1i_2\dots i_n$ in $C$
such that $i_1<i_2<\dots<i_d$ and
$i_{d+1}<\dots<i_{n-1}<i_n$. We will often identify such a minimal
representative $w$ with the vector $\bold i=(i_1,i_2,\dots,i_d)$ of its
first $d$ elements. The usual {\it Ehresmann--Bruhat order} on $S_n$
induces an order on the cosets of $S_n/(S_d\times S_{n-d})$. Given
two representatives, identified with $\bold i=(i_1,i_2,\dots,i_d)$ and
$\bold j=(j_1,j_2,\dots,j_d)$, respectively,
$\bold j$ is less or equal than $\bold i$
in this induced Bruhat order if and only if $j_\ell\le i_\ell$ for
all $\ell=1,2,\dots, d$.

Given a minimal representative $w$, we denote the corresponding
{\it Schubert variety} in the Gra\3mannian $\Gr_d(V)$ by $X(w)$. It is
well-known that $X(w)$ decomposes into the disjoint union
of {\it Schubert cells} which are indexed by elements
$\tau\in S_n/(S_d\times S_{n-d})$ with $\tau\le w$.
The multiplicity of a point $x$ in
$X(w)$ is constant on each Schubert cell.
Following \cite{\KrLaAA} we denote the
multiplicity of a point $x$ in the Schubert cell indexed by $\tau$
by $\mult_\tau X(w)$. In slight abuse of terminology
we will often call it the
``multiplicity of the point $\tau$ on the Schubert variety $X(w)$."

\smallskip
Let us now recall the multiplicity formula conjectured in \cite{\KrLaAA}.
We are given two elements $w$ and
$\tau$ in $S_n/(S_d\times S_{n-d})$.
In Conjecture~2 of \cite{\KrLaAA},
sets $S$ of reflections $s=(x,y)$, $1\le x\le d$, $d+1\le y\le
n$ (here we use standard transposition notation), are considered
with the property that
\roster
\item "(S1)" Any chain $s_1>s_2>\dots>s_t$ of commuting reflections,
all of them contained in $S$, satifies $w\ge\tau s_1\cdots s_t$ (in
the induced Bruhat order on $S_n/(S_d\times S_{n-d})$);
\item "(S2)" $S$ is maximal with respect to property (S1).
\endroster

Now we are in the position to formulate Conjecture~2 from
\cite{\KrLaAA}, which becomes a theorem by our proof in Section~3.

\proclaim{Theorem 1}
The multiplicity of the point $\tau$ on the Schubert variety $X(w)$
is given by
$$\mult _\tau X(w)=\vert\{S:S\text { satisfies (S1) and
(S2)}\}\vert.\tag2.1$$
\endproclaim

\midinsert
\vskip8pt
\vbox{
$$
\Gitter(10,22)(0,0)
%\Koordinatenachsen(-9,10)(0,0)
\PfadDicke{.5pt}
\Pfad(0,-1),2222222222222222222222\endPfad
\Pfad(-1,0),11111111111\endPfad
\PfadDicke{1.2pt}
\Pfad(1,9),22221222112212221111\endPfad
\Pfad(2,9),2212\endPfad
\Pfad(3,9),2\endPfad
\Pfad(4,9),222\endPfad
\Pfad(5,9),222222221221211\endPfad
\Pfad(6,9),22222221\endPfad
\Pfad(7,9),222222\endPfad
\Pfad(8,9),22222222212\endPfad
\Pfad(9,9),22222222\endPfad
%\SPfad(-1,9),11111111111\endSPfad
%\SPfad(1,10),112212221112211\endSPfad
\DickPunkt(1,9)
\DickPunkt(2,9)
\DickPunkt(3,9)
\DickPunkt(4,9)
\DickPunkt(5,9)
\DickPunkt(6,9)
\DickPunkt(7,9)
\DickPunkt(8,9)
\DickPunkt(9,9)
\DickPunkt(3,10)
\DickPunkt(3,12)
\DickPunkt(4,12)
\DickPunkt(7,15)
\DickPunkt(7,16)
\DickPunkt(9,17)
\DickPunkt(9,19)
\DickPunkt(9,20)
\DickPunkt(9,21)
\hskip4.5cm
$$
\centerline{\eightpoint Figure 1}
}
\vskip8pt
\endinsert

Next we recall the combinatorial interpretation of multiplicities in
terms of nonintersecting lattice paths from \cite{\KratBQ}, which is
a more or less straight-forward combinatorial translation of the
Rosenthal--Zelevinsky formula \cite{\RoZeAA} using the
Lindstr\"om--Gessel--Viennot theorem
\cite{\LindAA, Lemma~1}, \cite{\GeViAB, Theorem~1}.
As before, let $w$ and $\tau$ be two elements from $S_n/(S_d\times
S_{n-d})$, and identify them with $\bold i=(i_1,i_2,\dots,i_d)$,
$i_1<i_2<\dots<i_d$, and $\bold j=(j_1,j_2,\dots,j_d)$,
$j_1<j_2<\dots<j_d$, respectively.
Furthermore, we define the numbers $\kappa_q:=
\vert \{\ell:i_q<j_\ell\}\vert$. Then

\vskip-7pt
{\openup-1\jot
$$\align
\mult_\tau X(w)=&\#\text {\big(families $(P_1,P_2,\dots,P_d)$
of nonintersecting lattice paths,}\\
&\kern11pt
\text {where the
path $P_\ell$ runs from $(d+1-\ell,d)$
to }\\
&\kern11pt
\text {$(d-\kappa_{\sigma(\ell)},\kappa_{\sigma(\ell)}+i_{\sigma(\ell)}), \
\ell=1,2,\dots,d\big)$,}
\tag2.2
\endalign$$}
\vskip-7pt

\noindent
where $\sigma$ is some permutation in $S_d$.
See Figure~1. There, $d=9$, $\bold i=(4,6,7,13,14,17,\break 19,20,21)$
and $\bold j=(1,2,4,7,10,12,13,15,16)$. For this choice the vector of
the $\kappa_q$'s is $(6,6,5,2,2,0,0,0,0)$.
Figure~1 shows a typical family of paths as described in (2.2)
for this choice of $\bold i$ and
$\bold j$. The
permutation $\sigma$ is $674583129$ in this example.

At this point, there are two remarks to be made: First, in
\cite{\KratBQ}
the starting points of the paths are $(-\ell+1,\ell-1)$ and the end
points are $(-\kappa_{\ell},\kappa_{\ell}+i_{\ell})$, $\ell=1,2,\dots,d$ (the
latter in some order, determined by the permutation $\sigma$).
If we shift everything by $d$ units to the right then we obtain the
points $(d+1-\ell,\ell-1)$ and
$(d-\kappa_{\ell},\kappa_{\ell}+i_{\ell})$. Whereas now
the end points are in agreement, the starting points still differ slightly.
However, the arguments in Section~4
of \cite{\KratBQ} (and, in fact, figures such as Figure~3 in
\cite{\KratBQ}) show that
portions of paths below the horizontal line $y=d$ are forced and can
therefore be omitted. This means that we may replace the starting points
$(d+1-\ell,\ell-1)$ by the points $(d+1-\ell,d)$.
(In fact, Figure~1 shows exactly the result when the forced
portions of the paths in Figure~3 of \cite{\KratBQ} are cut off.)
Second, the order in which starting and end points are connected by
the nonintersecting paths is always the same, i.e., for fixed
$\bold i$ and $\bold j$ the permutation
$\sigma$ in (2.2) is always the same.

\subhead 3. Proof of the theorem\endsubhead
We will prove that the multiplicity formulas in Theorem~1 and (2.2)
are equivalent.

First we claim:
\proclaim{Claim 1}If a reflection $s=(x,y)$ is identified with the point
$(x,y)$ in the plane, then a set of reflections as described in {\rm(2.1)}
is the set of
all the lattice points with $y$-coordinates $>d$ on the paths of
a family of
paths as described in {\rm(2.2)}. In turn, given a
family of paths as described in {\rm(2.2)}, the set of lattice
points on the paths with $y$-coordinate $>d$ form a set of
reflections as described in {\rm(2.1)}, under the above identification of reflections
and points in the plane.
\endproclaim

To return to our example in Figure~1:
the set of lattice points on the paths with $y$-coordinate $>9$,
i.e., the set $\{(1,10),(1,11),(1,12),(1,13),(2,13),(2,14),\dots,
(9,21),\break(2,10),\dots,(3,12),\dots\}$, is a set of reflections with
the properties (S1) and (S2).

In fact, we are going to prove a more general claim. In order to be
able to formulate it, we have to explain the ``light-and-shadow
procedure with the sun in the south-east."
We will do this by considering an example.

Suppose that we are given a multiset $S$ of reflections, identified
with points in the plane as in Claim~1. For example, Figure~2.a
shows the multiset of reflections (points)
$$\multline
\{(2,13), (3,10), (3,10), (3,10), (3,11), (3,11), (4,10), (4,16), \\
(5,18), (5,18), (5,18),
(6,17), (6,18), (7,11), (7,16), (7,16), \\
(7,19), (8,21), (8,21), (8,21), (8,21), (9,13), (9,18) \}.
\endmultline$$

\midinsert
\vskip8pt
\vbox{
$$
\Gitter(10,22)(0,0)
%\Koordinatenachsen(-9,10)(0,0)
\PfadDicke{.5pt}
\Pfad(0,-1),2222222222222222222222\endPfad
\Pfad(-1,0),11111111111\endPfad
\DickPunkt(2,13)
\DickPunkt(4,16)
\DickPunkt(5,18)
\DickPunkt(3,11)
\DickPunkt(6,17)
\DickPunkt(7,19)
\DickPunkt(4,10)
\DickPunkt(7,11)
\DickPunkt(9,13)
\DickPunkt(6,18)
\DickPunkt(8,21)
\DickPunkt(3,10)
\DickPunkt(7,16)
\DickPunkt(9,18)
\Label\o{\eightpoint\ssize (2)}(3,11)
\Label\o{\eightpoint\ssize (4)}(8,21)
\Label\o{\eightpoint\ssize (3)}(5,18)
\Label\o{\eightpoint\ssize (3)}(3,10)
\Label\o{\eightpoint\ssize (2)}(7,16)
\hbox{\hskip7cm}
\Gitter(10,22)(0,0)
%\Koordinatenachsen(-9,10)(0,0)
\PfadDicke{.5pt}
\Pfad(0,-1),2222222222222222222222\endPfad
\Pfad(-1,0),11111111111\endPfad
\PfadDicke{1.2pt}
\Pfad(1,9),22221222112212221111\endPfad
\Pfad(2,9),2212\endPfad
\Pfad(3,9),2\endPfad
\Pfad(4,9),222\endPfad
\Pfad(5,9),222222221221211\endPfad
\Pfad(6,9),22222221\endPfad
\Pfad(7,9),222222\endPfad
\Pfad(8,9),22222222212\endPfad
\Pfad(9,9),22222222\endPfad
%\SPfad(-1,9),11111111111\endSPfad
%\SPfad(1,10),112212221112211\endSPfad
\Kreis(1,9)
\Kreis(2,9)
\Kreis(3,9)
\Kreis(4,9)
\Kreis(5,9)
\Kreis(6,9)
\Kreis(7,9)
\Kreis(8,9)
\Kreis(9,9)
\Kreis(3,10)
\Kreis(3,12)
\Kreis(4,12)
\Kreis(7,15)
\Kreis(7,16)
\Kreis(9,17)
\Kreis(9,19)
\Kreis(9,20)
\Kreis(9,21)
\DickPunkt(2,13)
\DickPunkt(4,16)
\DickPunkt(5,18)
\DickPunkt(3,11)
\DickPunkt(6,17)
\DickPunkt(7,19)
\DickPunkt(4,10)
\DickPunkt(7,11)
\DickPunkt(9,13)
\DickPunkt(6,18)
\DickPunkt(8,21)
\DickPunkt(3,10)
\DickPunkt(7,16)
\DickPunkt(9,18)
\Label\lo{\eightpoint\ssize (2)}(3,11)
\Label\o{\eightpoint\ssize (4)}(8,21)
\Label\lo{\eightpoint\ssize (3)}(5,18)
\Label\o{\eightpoint\ssize (3)}(3,10)
\Label\o{\eightpoint\ssize (2)}(7,16)
\Label\u{A_9}(1,9)
\Label\u{A_8}(2,9)
\Label\u{A_7}(3,9)
\Label\u{A_6}(4,9)
\Label\u{A_5}(5,9)
\Label\u{A_4}(6,9)
\Label\u{A_3}(7,9)
\Label\u{A_2}(8,9)
\Label\u{A_1}(9,9)
\Label\l{E_1}(3,10)
\Label\o{E_2}(3,12)
\Label\o{E_3}(4,12)
\Label\r{E_4}(7,15)
\Label\r{E_5}(7,16)
\Label\r{E_6}(9,17)
\Label\r{E_7}(9,19)
\Label\r{E_8}(9,20)
\Label\r{E_9}(9,21)
\hskip4.5cm
$$
\centerline{\eightpoint a. A multiset of points
\kern4cm
b. Light and shadow}
\centerline{\eightpoint Figure 2}
}
\vskip8pt
\endinsert

Next we suppose that there is a light source being located in the
bottom-right corner. The {\it shadow} of a point
$(x,y)$ is defined to be the set of points $(x',y')\in\Bbb R^2$ ($\Bbb
R$
denoting the set of real numbers) with $x'\le x$ and $y'\ge y$.
We consider the (bottom-right) {\it border} of the union of
the shadows of all the points of the multiset $S$.
We also include the shadows of the starting points $A_\ell=(d+1-\ell,d)$ and
the end points $E_\ell=(d-\kappa_\ell,\kappa_\ell+i_\ell)$,
$\ell=1,2,\dots,d$. This border is a
lattice path. We restrict our attention to the portion of this lattice
path between $A_1$ and $E_{\sigma(1)}$. (Here, as before, $\sigma$ is the
permutation as in (2.2) which describes how starting and end points are
connected in the case of nonintersecting lattice paths.
In our example in Figure~2,
$A_1$ and $E_{\sigma(1)}$ are the points $(9,9)$ and $(9,17)$, respectively.)
We remove all the points of the
multiset that lie on this path, including $A_1$ and $E_{\sigma(1)}$.
(In our example, we would remove $(9,9)$, $(9,13)$ and $(9,17)$.)
Then the
light and shadow procedure is repeated with the remaining
points. (That is, in the next step the roles of $A_1$ and
$E_{\sigma(1)}$ are played by $A_2$ and $E_{\sigma(2)}$, respectively,
etc.)
We stop after a total of $d$ iterations. (The result of applying this
procedure to the multiset in Figure~2.a is shown in Figure~2.b.)
It is obvious that at this point we will have obtained $d$
nonintersecting lattice paths, the $\ell$-th path connecting $A_\ell$
and $E_{\sigma(\ell)}$.

We are now ready to state:
\proclaim{Claim 2}If light-and-shadow with the sun in the south-east
is applied to a multiset $S$ of reflections satisfying {\rm(S1)} {\rm(}where
we again identify a reflection $s=(x,y)$ with the point
$(x,y)$ in the plane{\rm)}, then one obtains a family of
paths as described in {\rm(2.2)} which in addition cover all the
points of $S$. In turn, given a
family of paths as described in {\rm(2.2)}, any submultiset of the lattice
points on the paths with $y$-coordinate $>d$ forms a multiset of
reflections which satisfies {\rm(S1)}, under the above identification of reflections
and points in the plane.
\endproclaim

Claim~2 does indeed imply Claim~1: For suppose that we are given a set
$S$ of reflections (viewed as set of points) as described
in (2.1). Then the first
assertion of Claim~2 says that this set of points lies on a family of
paths as described in (2.2). Moreover, if $S$ were not
the complete set of lattice points on the paths with $y$-coordinate
$>d$, then we may add such a missing point, $(x,y)$ say, to $S$. The
second assertion of Claim~2 then says that $S\cup\{(x,y)\}$ is a set
of reflections satisfying (S1). Thus $S$ was not maximal, a
contradiction.
On the other hand, if we are given a family of paths as described in (2.2)
and consider the set $S$ of all lattice points on the paths with
$y$-coordinate $>d$, then the second assertion of Claim~2 says that
this is a set of reflections satisfying (S1). In addition, it is
maximal with respect to (S1). For suppose that it is not. Then we may
add another reflection, $(x,y)$ say, to $S$, thus obtaining $S'=S\cup
\{(x,y)\}$. Clearly, if we apply light-and-shadow to $S'$ then we will
not have exhausted all elements of $S'$ after these $d$
iterations (i.e., the $d$ paths obtained will not cover all elements
of $S'$). However, this is a contradiction to the first assertion of
Claim~2.

Claim~2 will be fully exploited in Section~4.

\smallskip
Let us call a set $\{(x_1,y_1),(x_2,y_2),\dots,(x_t,y_t)\}$ of points
with $x_1<x_2<\dots<x_t$ and $y_1>y_2>\dots>y_t$ a {\it chain}.
Furthermore, given a point $A=(a_1,a_2)$, let us define regions $R(A)$ by
$$R(A):=\{(x,y):x\le a_1,\ y>a_2\}.$$
(This is the region in the plane weakly to the left and strictly above
the point $A$.)

Claim~2 follows immediately from Claim~3 below. There, and in the
following, we assume tacitly that
any occurring set (multiset) of points is a subset (submultiset)
of the rectangle $\{(x,y):1\le x\le d,\ d+1\le y\le i_d\}$.

\proclaim{Claim 3}Let $\bold i$ and $\bold j$ be as before.
%Define the
%numbers $c_1,c_2,\dots,c_r$ (in decreasing order)
%to be the different values that are
%attained by the $\kappa_\ell$'s, i.e.,
%$$\{c_1,c_2,\dots,c_r\}=\{\kappa_1,\kappa_2,\dots,\kappa_d\}$$
%and $c_1>c_2>\dots>c_r$.
Then both the point multisets that satisfy {\rm(S1)}
and submultisets of lattice points with $y$-coordinate $>d$ taken
from a family of paths as described in {\rm(2.2)} can be characterized as
follows:
For any $q$ with $1\le q\le d$,
the maximal number of points that can be chosen from such a multiset
such that all of them are located inside $R(E_q)$ and in addition
form a chain
is at most $d-\kappa_q-q$. Here, as before, $E_q=(d-\kappa_q,\kappa_q+i_q)$.
\endproclaim

In the sequel, we will call the condition spelled out in the next-to-last
sentence of Claim~3 the {\it chain condition}.

Below we prove Claim~3, which in fact means to prove four assertions,
labelled A1--A4.

\smallskip
A1. {\it Any submultiset of lattice points with $y$-coordinate $>d$ taken
from a family of paths as described in {\rm(2.2)} satisfies the chain
condition}. This is obvious once one observes that $d-\kappa_q-q$ is the
number of lattice paths in the family that start strictly to the left of $E_q$
and terminate weakly to the right of $E_q$ (and, enforcedly, pass above
$E_q$; cf\. Figure~1).

\smallskip
A2. {\it If a multiset of lattice points satisfies the chain condition
then it is a submultiset of the lattice points with $y$-coordinate
$>d$ of a family of paths as described in} (2.2). This is also more or
less ``obvious.'' The only matter is notation. Probably the most
convenient way to prove this rigorously is by induction on $d$.

\midinsert
\vskip8pt
\vbox{
$$
\Gitter(10,22)(0,0)
%\Koordinatenachsen(-9,10)(0,0)
\PfadDicke{.5pt}
\Pfad(0,-1),2222222222222222222222\endPfad
\Pfad(-1,0),11111111111\endPfad
\PfadDicke{1.2pt}
\Pfad(8,21),1\endPfad
\Pfad(8,20),1\endPfad
\Pfad(8,18),12\endPfad
\Pfad(9,9),22222222\endPfad
%\SPfad(-1,9),11111111111\endSPfad
%\SPfad(1,10),112212221112211\endSPfad
\Kreis(1,9)
\Kreis(2,9)
\Kreis(3,9)
\Kreis(4,9)
\Kreis(5,9)
\Kreis(6,9)
\Kreis(7,9)
\Kreis(8,9)
\Kreis(9,9)
\Kreis(3,10)
\Kreis(3,12)
\Kreis(4,12)
\Kreis(7,15)
\Kreis(7,16)
\Kreis(9,17)
\Kreis(9,19)
\Kreis(9,20)
\Kreis(9,21)
\DickPunkt(2,13)
\DickPunkt(4,16)
\DickPunkt(5,18)
\DickPunkt(3,11)
\DickPunkt(6,17)
\DickPunkt(7,19)
\DickPunkt(4,10)
\DickPunkt(7,11)
\DickPunkt(9,13)
\DickPunkt(6,18)
\DickPunkt(8,21)
\DickPunkt(3,10)
\DickPunkt(7,16)
\DickPunkt(9,18)
\Label\lo{\eightpoint\ssize (2)}(3,11)
\Label\o{\eightpoint\ssize (4)}(8,21)
\Label\lo{\eightpoint\ssize (3)}(5,18)
\Label\o{\eightpoint\ssize (3)}(3,10)
\Label\o{\eightpoint\ssize (2)}(7,16)
\Label\u{A_9}(1,9)
\Label\u{A_8}(2,9)
\Label\u{A_7}(3,9)
\Label\u{A_6}(4,9)
\Label\u{A_5}(5,9)
\Label\u{A_4}(6,9)
\Label\u{A_3}(7,9)
\Label\u{A_2}(8,9)
\Label\u{A_1}(9,9)
\Label\l{E_1}(3,10)
\Label\o{E_2}(3,12)
\Label\o{E_3}(4,12)
\Label\r{E_4}(7,15)
\Label\r{E_5}(7,16)
\Label\r{E_6}(9,17)
\Label\r{E_7}(9,19)
\Label\r{E_8}(9,20)
\Label\r{E_9}(9,21)
\Label\l{E'_7}(8,18)
\Label\l{E'_8}(8,20)
\Label\l{E'_9}(8,21)
\hskip4.5cm
$$
\centerline{\eightpoint Figure 3}
}
\vskip8pt
\endinsert

For $d=1$ the assertion is obvious (the quantity $d-\kappa_q-q$ being 0 for
$d=q=1$). Let us now assume that we have
already proved the assertion for $d$. Given $\bold
i=(i_1,i_2,\dots,i_{d+1})$ and $\bold j=(j_1,j_2,\dots,j_{d+1})$ and
a multiset of lattice points satisfying the chain condition, we
apply light-and-shadow (with respect to the starting and end points
determined by $\bold i$ and $\bold j$). We restrict our attention to
the rightmost strip of the picture, i.e., the region of points with
$x$-coordinate between $d$ and $d+1$, see Figure~3. There, we have
chosen $d=8$, $\bold i=(4,6,7,13,14,17,19,20,21)$
and $\bold j=(1,2,4,7,10,12,13,15,16)$. The starting and end points
determined by $\bold i$ and $\bold j$ are indicated by circles. The
multiset of points is indicated by bold dots, multiplicities being
indicated by the numbers in parentheses. (This is in fact the same
example as in Figure~2. The path pieces should be ignored for the
moment.)

Let $k$ be minimal such that $i_k\ge j_{d+1}$. (In our example we have
$k=6$.) Then the end points with
$x$-coordinate $d+1$ are $E_k,E_{k+1},\dots,E_{d+1}$. Clearly, under
light-and-shadow, $E_k$ is connected with $A_1$. The path portions
leading to the other end points $E_{k+1},\dots,E_{d+1}$ hit the
vertical line $y=d$ the last time in the points
$E'_{k+1},\dots,E'_{d+1}$, say (see Figure~3). It is easy to see that for any
$q$ with $k+1\le q\le d+1$
the maximal number of points that can be chosen from such a multiset
such that they form a chain and all of them are located inside $R(E'_q)$
is at most $d+1-q$. Now we apply the induction hypothesis to
$\bold i'=(i_1,\dots,i_{k-1},i'_{k+1}-1,\dots,i'_{d+1}-1)$ and $\bold
j'=(j_1,j_2,\dots,j_d)$,
where $i'_{\ell}$ denotes the $y$-coordinate of $E'_{\ell}$,
$\ell=k+1,\dots,d+1$.
It should be observed that, up to a vertical shift of 1 unit,
the starting points determined by $\bold i'$
and $\bold j'$ are $A_2,A_3,\dots,A_{d+1}$, whereas the corresponding
end points are $E_1,\dots,E_{k-1},E'_{k+1},\dots,E'_{d+1}$. By the
above consideration, the multiset satisfies the chain condition with
respect to these new starting and end points. The
induction hypothesis then guarantees that light-and-shadow yields a
family of paths connecting the (new) starting points with the (new) end
points, thereby covering all (remaining)
elements of the multiset. This family
of paths is finally concatenated with the path portions that we
already obtained in the strip between the vertical lines $y=d$ and
$y=d+1$.

\smallskip
A3. {\it Any multiset of reflections satisfying
{\rm(S1)} satisfies the chain condition}.
Suppose we are given a multiset of reflections which satisfies (S1)
but does not satisfy the chain condition.
Then for some $q$ there is a chain of $d-\kappa_q-q+1$ reflections
from the multiset which, when viewed as points in the plane, are all
located inside $R(E_q)$.

Let the reflections in the chain be $s_1,s_2,\dots,s_{d-\kappa_q-q+1}$.
Let us consider a reflection in the chain, $(x,y)$ say.
Since $(x,y)\in R(E_q)$ we have $x\le d-\kappa_q$. Furthermore, we have
$$(\tau s_1\cdots s_{d-\kappa_q-q+1})(x)=\tau(y),$$
where as before $\tau$ is identified with $\bold j$, i.e.,
$\tau=j_1j_2\dots j_n$ with $j_1<j_2<\dots<j_d$ and
$j_{d+1}<\dots<j_{n-1}<j_n$.
Since $(x,y)$ is contained in $R(E_q)$ we have $y>\kappa_q+i_q\ge d$.
Because of $j_{d+1}<\dots<j_{n-1}<j_n$, this implies $\tau(y)\ge
\tau(\kappa_q+i_q+1)$.

We claim that $\tau(\kappa_q+i_q+1)=i_q+1$. This is seen as follows.
Taking into account the trivial fact that the set of values
$\{j_{d+1},j_{d+2},\dots,j_{n}\}$ is equal to the complement of
$\{j_{1},j_{2},\dots,j_{d}\}$ in $\{1,2,\dots,n\}$, a
value $\tau(y)=j_y$ for $y>d$ is {\it characterized\/} by
$$j_y=(y-d)+\vert\{\ell:j_{\ell}<j_{y}\}\vert.$$
Thus we may verify our claim by setting $y=\kappa_q+i_q+1$ and substituting
$i_q+1$ for $j_{\kappa_q+i_q+1}$ in this equation. Indeed, we have
$i_1+1=(\kappa_q+i_q+1-d)+(d-\kappa_q)$. Hence, we have $\tau(y)>i_q$.

In summary, we have found $d-\kappa_q-q+1$ values of $x$ such that
$$(\tau s_1\cdots s_{d-\kappa_q-q+1})(x)>i_q,\tag3.1$$
all of which are $\le d-\kappa_q$.
Moreover, if $d-\kappa_q<x\le d$, then we have
$$(\tau s_1\cdots s_{d-\kappa_q-q+1})(x)=\tau(x)=j_x>i_q.$$
Hence, in total we found $(d-\kappa_q-q+1)+\kappa_q=d-q+1$ values
$x$ for which (3.1) holds. If we recall that we also always identify $w$
with $\bold i$, i.e., $w=i_1i_2\dots i_n$, then
this is a contradiction to (S1).

\smallskip
A4. {\it If a multiset of lattice points satisfies the chain condition
then, if viewed as a multiset of reflections, it also satisfies}
(S1).
Consider a chain of $t$ points of the multiset, and view them as
reflections $s_1>s_2>\dots>s_t$.
Still identifying $w$ and $\bold i$, we observe that
the inequality $w\ge\tau s_1\cdots s_t$ is equivalent to the
inequality
$$\big\vert\{x:1\le x\le d\text{ and }(\tau s_1\cdots
s_t)(x)>i_q\}\big\vert \le d-q\tag3.2$$
to hold for $1\le q\le d$.

A careful examination of the arguments in A3 shows that they actually
prove
$$\big\vert\{x:1\le x\le d\text{ and }(\tau s_1\cdots
s_t)(x)>i_q\}\big\vert =\big\vert\{\ell:s_\ell\in R(E_q)\}\big\vert +
\kappa_q.$$
By assumption, our multiset of points satisfies the chain condition,
hence $\big\vert\{\ell:s_\ell\in R(E_q)\}\big\vert \le
d-\kappa_q-q$. Clearly, this implies (3.2), as desired.

\subhead 4. A formula for the Hilbert series of the tangent cone at
a point\endsubhead
Now the full significance of Claim~2 can be revealed. Briefly, it allows the
formulation of a version of Conjecture~1 in \cite{\KrLaAA}
which has the
advantage of being {\it efficient}, as it reduces the computation of the
Hilbert function to a {\it finite} problem. More precisely, we can express
the {\it Hilbert series} in form of a finite summation. This form of
the conjecture is the analogue of, say, formulas for the Hilbert
series as in \cite{\CoHeAA, Eq.~(1)}, \cite{\KrPrAA, p.~1021, line~11}
or in \cite{\KrRuAA, Theorem~1}.

In order to formulate this equivalent form, we need to introduce some
notation.
A point in a lattice path $P$ which is the end point of a horizontal
step and at the same time the starting point of a vertical step will
be called an {\it east-north turn\/} ({\it EN-turn\/} for short) of the
lattice path $P$. For example, the EN-turns of the leftmost lattice paths in
Figure~1 are $(2,13)$, $(4,16)$ and $(5,18)$.
We write $\EN(P)$ for the number of
NE-turns of $P$. Also, given a family $\bold P=(P_1,P_2,\dots,P_d)$ of
paths $P_\ell$, we write $\EN(\bold P)$ for the number $\sum _{\ell=1}
^{d}\EN(P_\ell)$ of all EN-turns in the family.
By $\Bbb P^+_L(\bold
A\to \bold E)$ we denote the set of all families $(P_1,P_2,\dots,P_d)$
of nonintersecting lattice paths, where $P_\ell$ runs from $A_\ell$ to
$E_\ell$.
Finally, given any weight function $\mu$ defined on a set $\Cal M$, by
the generating function $\GF(\Cal M;\mu)$ we mean $\sum _{x\in\Cal M}
^{}\mu(x)$.

\proclaim{Theorem 2} Conjecture~1 from \cite{\KrLaAA}
is equivalent to saying that the
Hilbert series of the tangent cone to $X(w)$ at $\tau$ is equal to
$$\frac {\GF(\Bbb P^+(\bold A\to \bold E);z^{\EN(.)})}
{(1-z)^{\sum _{\ell=1} ^{d}i_\ell-\binom {d+1}2}},
\tag4.1$$
with $A_\ell=(d+1-\ell,d)$ and
$E_\ell=(d-\kappa_\ell,\kappa_\ell+i_\ell)$,
$\ell=1,2,\dots,d$, as before.
\endproclaim

\demo{Proof} We can simply copy the corresponding proof in
\cite{\KrPrAA, first proof of Theorem~2}.

According to the conjecture, the dimension
of the $m$-th homogeneous component of the tangent cone is equal to
the number of multisets of cardinality $m$ which satisfy (S1).
Following \cite{\KrLaAA} we denote this dimension by $h_{TC_\tau
X(w)}(m)$.

Let $S$ be such a multiset. We apply light-and-shadow to it. By
Claim~2, we obtain a family $(P_1,P_2,\dots,P_d)$ of paths as described in (2.2).
Each path $P_\ell$ contains a few (possibly multiple) points of $S$.
However, in each EN-turn of $P_\ell$ there has to be {\it at least
one\/} element of $S$, $\ell=1,2,\dots,d$.
Therefore, and because of the second assertion of Claim~2,
given such a family $(P_1,P_2,\dots,P_d)$ of paths as described in
(2.2) with a total number of
exactly $t$ EN-turns, there are exactly $\binom{T+m-t-1}{m-t}$
multisets $S$ of cardinality $m$ that reduce to
$(P_1,P_2,\dots,P_d)$ under light and shadow, where
$$T=\sum _{\ell=1} ^{d}((d-\kappa_\ell)+(\kappa_\ell+i_\ell))-
\sum _{\ell=1} ^{d}((d+1-\ell)+d)=
\sum _{\ell=1} ^{d}i_\ell-\binom {d+1}2$$
is the total number of lattice points with $y$-coordinate $>d$
on the lattice paths
$P_1,P_2,\dots,P_d$. (It is independent of the path family.)

Hence, if we let $h_t$ denote the number of all
families $(P_1,P_2,\dots,P_d)$ of
paths as described in (2.2) with a total number of
exactly $t$ EN-turns, we obtain for the Hilbert series,
$$\align \sum _{m=0} ^{\infty}h_{TC_\tau X(w)}(m)\, z^m&=
\sum _{m=0} ^{\infty}\bigg(\sum _{t=0}
^{m}\binom{T+m-t-1}{m-t}h_t\bigg)z^m\\
&=\sum _{t=0} ^{\infty}h_t\sum _{m=t}
^{\infty}\binom{T+m-t-1}{m-t}z^m\\
&=\sum _{t=0} ^{\infty}h_tz^t\sum _{m=0}
^{\infty}\binom{T+m-1}{m}z^m\\
&=\frac {\sum _{t=0} ^{\infty}h_tz^t} {(1-z)^T}.
\endalign$$
Now, the generating function $\sum _{t=0} ^{\infty}h_tz^t$ is exactly
the numerator in (4.1).
This proves the theorem.\quad \quad \qed
\enddemo

\remark{Remarks} (1) If true, formula (4.1) implies the
Rosenthal--Zelevinsky formula. For, the multiplicity $\mult_\tau X(w)$
is equal to the {\it numerator} of the
Hilbert series of the tangent cone to $X(w)$ at $\tau$,
{\it evaluated at $z=1$}. But, by (4.1), this is
exactly the {\it number} of all families of nonintersecting lattice
paths in $\Bbb P^+(\bold A\to \bold E)$, i.e., of all path families as
described in (2.2). As we already remarked earlier, the combinatorial
interpretation of the multiplicity in terms of nonintersecting
lattice paths as given in (2.2) is equivalent to the
Rosenthal--Zelevinsky formula.

(2)  Unfortunately, all the results that
have been found so far on the enumeration of nonintersecting lattice
paths with respect to turns (see \cite{\KratBE, \KrPrAA, \KrRuAA,
\KulkAC, \ModaAA})
do not cover the above case, because the
location of the starting and end points is quite unusual. This means
that, up to now, there is no compact formula (a determinant, or
whatever) for
$\GF(\Bbb P^+(\bold A\to \bold E);z^{\EN(.)})$.
\endremark

\enddocument